\documentstyle[11pt, epsfig]{article}
\input amssym.def
\title{Solutions of the Einstein-Dirac Equation on Riemannian 3-Manifolds with Constant Scalar Curvature.  
\footnote{Supported by the SFB 288 of the DFG.}}
\author{Thomas Friedrich}
\date{{\small \it Humboldt-Universit\"at zu Berlin, Institut f\"ur Reine Mathematik,\\
Ziegelstra\ss e 13a, D-10099 Berlin, Germany}\\}

\begin{document}

\newcommand{\D}{\displaystyle}
\newcommand{\upsp}{\phantom{l}}
\newcommand{\downsp}{\phantom{q}}

\maketitle

\mbox{} \hrulefill \mbox{}\\

\newcommand{\vol}{\mbox{vol} \, }
\newcommand{\grad}{\mbox{grad} \, }

\begin{abstract}  \noindent This paper contains a classification of all 3-dimensional
manifolds with constant scalar curvature $S \not= 0$ that carry a non-trivial
solution of the Einstein-Dirac equation.
\end{abstract}

{\small
{\it Subj. Class.:} Differential Geometry.
 }\\

\setcounter{section}{0}

\mbox{} \hrulefill \mbox{}\\

\section{Introduction}

Consider a Riemannian spin manifold of dimension $n \ge 3$ and denote by $D$ the Dirac operator 
acting on spinor fields. A solution of the Einstein-Dirac equation is a spinor field $\psi$ solving the equations

$$  Ric - \frac{1}{2} \, \, S \cdot g  =  \pm \frac{1}{4}  \, \, T_{\psi} \quad , \quad D(\psi) = \lambda \psi . $$

Here $S$ denotes the scalar curvature of the space, $\lambda$ is a real constant and $T_{\psi}$ is the 
energy-momentum tensor of the spinor field $\psi$ defined by the formula
$$ T_{\psi} (X,Y) =(X \cdot \nabla_Y \psi + Y \cdot \nabla_X \psi , \psi ) . $$

The scalar curvature $S$ is related to the eigenvalue $\lambda$ and the length of the spinor field
$\psi$ by the formula
$$ S= \pm \frac{\lambda}{n-2} \vert \psi \vert^2 . $$

In [KimF] we introduced the weak Killing equation for a spinor field $\psi^*$:
$$ \nabla_X \psi^* = \frac{n}{2(n-1)} \, \, dS (X) \psi^* + \frac{2 \lambda}{(n-2)S} \, \, 
Ric (X) \cdot \psi^* - \frac{\lambda}{n-2} \, \, X \cdot \psi^* + \frac{1}{2(n-1)S} \, \, X \cdot dS \cdot \psi^* $$

Any weak Killing spinor $\psi^*$ (WK-spinor) yields a solution $\psi$ of the Einstein-Dirac 
equation after normalization
$$\psi = \sqrt{\frac{(n-2)\vert S \vert}{\vert \lambda \vert \vert \psi^* \vert^2}} \, \, \psi^* . $$

In fact, in dimension $n=3$ the Einstein-Dirac equation is essentially equivalent to the weak Killing equation 
(see [KimF]). Up to now the following 3-dimensional Riemannian manifolds admitting WK-spinors are known:
\begin{enumerate}
\item the flat torus $T^3$ with a parallel spinor;
\item the sphere $S^3$ with a Killing spinor;
\item two non-Einstein Sasakian metrics on the sphere $S^3$ admitting WK-spinors.
The scalar curvature of these two left-invariant metrics equals $S= 1 \pm \sqrt{5}$.
\end{enumerate}

The aim of this paper is to classify all Riemannian 3-manifolds with constant scalar curvature
and admitting a solution of the Einstein-Dirac equation. In particular, we will prove
 the existence of a one-parameter family of left-invariant metrics on 
$S^3$ with WK-spinors. This family contains the two non-Einstein Sasakian 
metrics with WK-spinors on $S^3$, 
but does not contain the standard sphere $S^3$ with Killing spinors. Moreover, any simply-connected, complete 
Riemannian manifold $N^3 \not= S^3$ with WK-spinors and constant scalar curvature is isometric to a space 
of this one-parameter family. In order to formulate the result precisely, 
we fix real parameters $K,L,M \in {\Bbb R}$ and denote by $N^3 (K,L,M)$ 
the 3-dimensional, simply-connected and oriented Riemannian manifold defined by the following structure equations:
$$ \omega_{12} = K \sigma^3 \quad , \quad \omega_{13}= L \sigma^2 \quad , \quad \omega_{23} = M \sigma^1  , $$

or, equivalently:

$$ d \sigma^1 =(L-K) \sigma^2 \wedge \sigma^3 \quad , \quad d \sigma^2 = (M+K) \sigma^1 \wedge \sigma^3 \quad , 
\quad d \sigma^3 = (L-M) \sigma^1 \wedge \sigma^2 . $$

The 1-forms $\sigma^1, \sigma^2 , \sigma^3$ are the dual forms of an orthonormal frame of vector fields. 
Using this frame the Ricci tensor of $N^3 (K,L,M)$ is given by the matrix
$$ Ric = \left( \begin{array}{ccc}
-2KL & 0 & 0\\
0 & 2KM & 0 \\
0 & 0 & -2LM
\end{array} \right) . $$

{\bf Main Theorem:} {\it Let $N^3 \not= S^3$ be a complete, simply-connected Riemannian manifold with 
constant scalar curvature $S \not= 0$.  If $N^3$ admits a WK-spinor, then $N^3$ is isometric to 
$N^3 (K, L,M)$ and the parameters are a solution of the equation 
$$ - K^2 L(L-M)^2 M+ L^3 M^3 + KL^2 M^2 (M-L) + K^3 (L-M)(L+M)^2=0 \hspace{1cm} (*)$$

Conversely, any space $N^3 (K,L,M)$ such that  $(K,L,M) \not= (0,0,0)$ is a solution of  $(*)$ admits two 
WK-spinors for one and only one WK-number $\lambda$. With respect to the fixed orientation of $N^3 (K,L,M)$ 
we have the two cases:
$$ \lambda = + \frac{S}{2 \sqrt{2}} \, \, \sqrt{ \frac{S}{S^2 - 2|Ric|^2}} \quad \quad \mbox{if \,  $-K <M$}$$
$$ \lambda = - \frac{S}{2 \sqrt{2}} \, \, \sqrt{ \frac{S}{S^2 - 2|Ric|^2}} \quad \quad  \mbox{if \,  $M < -K$} . $$

\bigskip

The spaces $N^3 (K,L,M)$ are isometric to $S^3$ equipped with a left-invariant metric.\\}

\bigskip

{\bf Remark:} If the parameters $K=M$ coincide, the solution of the equation $(*)$ is given by 
$$ L = \frac{1}{4} K(1 - \sqrt{5}) \quad , \quad L= \frac{1}{4} K(1+ \sqrt{5})  $$

and we obtain the Ricci tensors

$$ Ric = \left( \begin{array}{ccc}
\frac{1}{2} K^2 (\sqrt{5}-1) & 0 & 0\\
0 & 2K^2 & 0\\
0 & 0 & \frac{1}{2} K^2 (\sqrt{5}-1) \end{array} \right)$$

or

$$ Ric = \left( \begin{array}{ccc}
- \frac{1}{2} K^2 (1+ \sqrt{5}) & 0 & 0\\
0 & 2K^2 & 0\\
0 & 0 & - \frac{1}{2} K^2 (1+ \sqrt{5}) \end{array} \right) . $$

\bigskip

The non-Einstein-Sasakian metrics on $S^3$ occur for the parameter $K=1$ (see [KimF]).\\

\newfont{\graf}{eufm10}
\newcommand{\sod}{\mbox{\graf so}(3)}

{\bf Remark:} Using the standard basis of the Lie algebra $\sod$ we can write the left-invariant metric 
of the space $N^3 (K,L,M)$ in the following way:
$$ \left( \begin{array}{ccc}
\frac{1}{|M-L||K+M|} & 0 & 0\\
\\
0 & \frac{1}{|K-L||M-L|} & 0\\
\\
0 & 0 & \frac{1}{|K-L||K+M|} \end{array} \right) . $$

The equation $(*)$ is a homogeneous equation of order six. The transformation $(K,L,M) \to 
(\mu K, \mu L, \mu M)$ corresponds to a homothety of the metric. Therefore - up to a homothety - 
the moduli space of solutions is a subset of the real projective space ${\Bbb P}^2 ({\Bbb R})$ given 
by the equation $(*)$. This subset is a configuration of six curves in ${\Bbb P}^2 ({\Bbb R})$ connecting the three points $[K:L:M] = [1:0:0], [0:1:0], [0:0:1]$ corresponding to flat metrics. 

\begin{center}

\[
\epsfig{figure=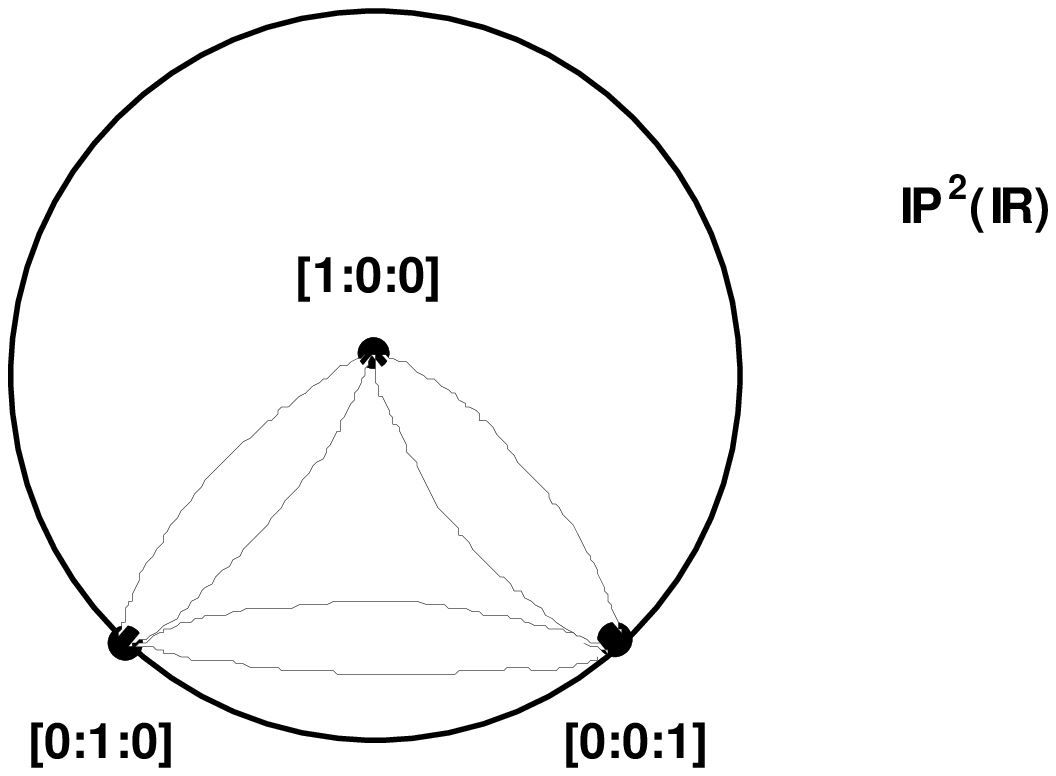,width=8cm}
\]

\end{center}

In particular, we have constructed two paths of solutions of the Einstein-Dirac equation deforming 
the non-Einstein Sasakian metrics on $S^3$.\\

\section{The integrability condition for the Einstein-Dirac equation in dimension $n=3$.}

The spinor bundle of a 3-dimensional Riemannian manifold is a complex vector bundle of dimension
two. Moreover, there exists a quaternionic structure commuting with the Clifford multiplication
by real vectors (see [F]). Consequently, in case of a real WK-number $\lambda$, the corresponding
space of WK-spinors is a quaternionic vector space. In the spinor bundle let us introduce
the metric connection $\nabla^{\lambda}$ given by the formula 
$$  \nabla^{\lambda}_X \psi := \nabla_X \psi - \frac{3}{4} dS(X) \psi - 
\lambda \Big\{ \frac{2}{S} Ric (X) - X \Big\} \cdot \psi - \frac{1}{4S} X \cdot dS \cdot \psi $$

and denote by $\Omega^{\lambda}$ its curvature form. Then we obtain the following\\

{\bf Proposition 1:} {\it Let $N^3$ be a simply-connected 3-dimensional Riemannian
manifold and suppose that the scalar curvature $S \not= 0$ does not vanish.
Then the following conditions are equivalent:
\begin{enumerate}
\item $N^3$ is a non-trivial solution of the Einstein-Dirac equation with real
eigenvalue $\lambda$;
\item $N^3$ admits a WK-spinor with real WK-number $\lambda$;
\item $N^3$ admits two WK-spinors with real WK-number $\lambda$;
\item The curvature $\Omega^{\lambda} \equiv 0$ vanishes identically.
\end{enumerate}}

If the scalar curvature $S \not= 0$ is constant, the condition $\Omega^{\lambda} \equiv 0$
 has been investigated and yields algebraic equations involving the
Ricci tensor and its covariant derivative (see [KimF], Theorem 8.3.). 
In order to formulate the integrability condition, we denote by $X \times Y$
the vector cross product of two vectors $X,Y \in T(N^3)$. For brevity,
let us introduce the endomorphism $T:T (N^3) \to T(N^3)$ given by the 
formula
$$ T(X) = \sum\limits^3_{i=1} e_i \times (\nabla_{e_i} Ric) (X) , $$

which will be used in the proof of the main Theorem.\\

{\bf Theorem 1 (see [KimF]):} {\it Let $N^3$ be a simply-connected 3-dimensional
Riemannian manifold with constant scalar curvature $S \not= 0$. 
$N^3$ admits a solution of the Einstein-Dirac equation with real eigenvalue
$\lambda$ if and only if the following three conditions are satisfied:}

\begin{enumerate}
\item $8 \lambda^2 \{ S^2 - 2 |Ric |^2 \} = S^3 ; $
\item $8 \lambda^2 \{ S Ric (X) - 2 Ric \circ Ric (X) \} - 4 \lambda ST(X)- S^2 Ric (X) =0;$
\item $8 \lambda^2 \{ 2 Ric (X) - SX \} \times \{ 2 Ric (Y) - SY \} + 8 \lambda S \{ (\nabla_X Ric )(Y) - 
(\nabla_Y Ric)(X) \}  $\\ [0.5em]
$+ S^3 X \times Y =  2 S^2 \sum\limits_{i<j} \{ R_{jY} \delta_{iX} + R_{iX} \delta_{jY} \} e_i \times {e_j}.$
\end{enumerate}

\section{Proof of the Main Theorem}

We fix an orthonormal frame $e_1, e_2, e_3$ of vector fields on $N^3$ consisting of eigenvectors
of the Ricci tensor:
$$ Ric = \left( \begin{array}{ccc} A&0&0\\0&B&0\\0&0&C \end{array} \right) . $$

Denote by $\sigma^1, \sigma^2, \sigma^3$ the dual frame and consider the connection forms 
$\omega_{ij}= \langle \nabla e_i , e_j \rangle$ of the Levi-Civita connection. The structure
equations of the Riemannian manifold $N^3$ are
$$ d \omega_{12} = \omega_{13} \wedge \omega_{32} +\frac{C-A-B}{2}  \, \, \sigma^1 \wedge \sigma^2$$
$$ d \omega_{13} = \omega_{12} \wedge \omega_{23} +\frac{B-A-C}{2}  \, \, \sigma^1 \wedge \sigma^3$$
$$ d \omega_{23} = \omega_{21} \wedge \omega_{13} +\frac{A-B-C}{2}  \, \, \sigma^2 \wedge \sigma^3$$

and the covariant derivative $\nabla Ric$ is given by the matrix of 1-forms
$$ \nabla Ric = \left( \begin{array}{ccc} \smallskip dA & (A-B) \omega_{12} & (A-C) \omega_{13}\\ \smallskip
(A-B) \omega_{12} & dB & (B-C) \omega_{23}\\
(A-C) \omega_{13} & (B-C) \omega_{23} & dC \end{array} \right) . $$

Using the third equation of Theorem 1 we obtain
$$ \langle (\nabla_{e_i} Ric)(e_j) - (\nabla_{e_j} Ric)(e_i), e_i \rangle =0 $$

and, consequently,
$$ dA(e_1)=dA(e_2)=dB(e_1)=dB(e_3)=dC(e_2)=dC (e_3)=0 . $$

Since $A+B+C=S$ is constant, we conclude that any eigenvalue $A,B,C$ of the Ricci tensor is constant, too. 
The second equation of Theorem 1 yields the condition that all elements outside the diagonal of
the (1,1)-tensor $T$ are zero:
$$ (A-B) \omega_{12} (e_1) =0= (A-B) \omega_{12} (e_2)$$
$$ (C-A) \omega_{13} (e_1) =0= (C-A) \omega_{13} (e_3)$$
$$ (B-C) \omega_{23} (e_2) =0= (B-C) \omega_{23} (e_3) . $$

First, we discuss the generic case that $A,B,C$ are pairwise different.
Then there exist numbers $K,L,M$ such that
$$\omega_{12} = K \sigma^3 \quad , \quad \omega_{13} = L \sigma^2 \quad , 
\quad \omega_{23} =M \sigma^1 . $$

The parameter triples $\{A,B,C \}$ and $\{K,L,M \}$ are related via the structure equations
by the formulas
$$ A=-2KL \quad , \quad B=2KM  \quad , \quad C= - 2 LM . $$

The first and second equation of Theorem 1 become equivalent to the following
system of algebraic equations:

\begin{itemize}
\item[1'.] $\displaystyle \lambda = \pm \frac{S}{2\sqrt{2}} \sqrt{\frac{S}{S^2 - 2 \vert Ric \vert^2}}$ ;
\item[2'.] \mbox{} \, $2S (S^2  - 2|Ric|^2) \{ (A-C) L+(B-A)K \}^2 \, \, =\, \, S(SA-2A^2)-A(S^2-2|Ric|^2)$
\begin{eqnarray*}2S (S^2 - 2|Ric|^2) \{ (C-B)M +(A-B)K \}^2 &=&S(SB-2B^2)-B(S^2-2|Ric|^2)\\ [0.4em]
2S (S^2 - 2|Ric|^2) \{ (B-C)M +(C-A)L \}^2 &=&S(SC-2C^2)-C(S^2-2|Ric|^2) . 
\end{eqnarray*}
\end{itemize}

We solve this system of algebraic equations with respect to the parameters
$K,L,M$. It turns out that the equations 2'. can be written in the form
$$ P_i (K,L,M) \cdot Q(K,L,M) =0 , $$

$(1 \le i \le 3)$, where the polynomials $P_1, P_2, P_3$ and $Q$ are given by the 
formulas
\begin{eqnarray*}
P_1 (K,L,M) &=& (-KL^2 + L^2M+K^2 (L+M))^2\\
P_2 (K,L,M) &=& (KM^2 + LM^2 + K^2 (L+M))^2\\
P_3 (K,L,M) &=& (LM(-L+M) +K(L^2+M^2))^2\\
Q  \hspace{0.1cm}  (K,L,M) &=& -K^2L(L-M)^2 M+L^3M^3+KL^2M^2(M-L) + K^3(L-M)(L+M)^2 .
\end{eqnarray*}

The real solutions of $P_1= P_2= P_3 =0$ are the pairs $\{K=0, L=0\}$ (the flat metric) and 
$\{K=M,L=-M\}$ (the space of positive constant curvature). Therefore, we 
proved  that a 3-dimensional complete, simply-connected manifold $N^3$ with
constant scalar curvature $S \not= 0$ and different eigenvalues of the Ricci tensor is 
isometric to one of the spaces $N^3 (K,L,M)$, where the parameters $K,L,M$ are
solutions of the equation $Q(K,L,M)=0$. These spaces satisfy the conditions
1. and 2. of Theorem 1 and, moreover, a simple computation yields the result
that condition 3. of Theorem 1 is satisfied, too. We next discuss the case 
that two of the eigenvalues $A,B,C$ coincide, for
example, $A=C \not= B$. Then we obtain again
$$ \omega_{12} = K \sigma^3 \quad , \quad \omega_{23} =M \sigma^1 , $$

but there is no condition for the connection form $\omega_{13}$. We compute
the matrix of the $(1,1)$-tensor $T$:
$$ T= \left( \begin{array}{ccc}
(B-C)K & 0&0\\ [0.2em]
0& (C-B)(K+M) & 0\\ [0.2em]
0&0 & (B-C)M
\end{array} \right) . $$

Since the scalar curvature $S$ as well as the eigenvalues $A=C,B$ of the
Ricci tensor are constant, the second equation of Theorem 1 yields
that $K$ and $M$ are constant and, moreover, coincide:
$$K=M= \, \mbox{const.}  $$

In case of $K=M=0$ we have $\omega_{12}= \omega_{23}=0$ and
$A=C$. In particular, the Ricci tensor is parallel, $\nabla Ric =0$. Therefore,
in this case $N^3$ is a Ricci-parallel 3-dimensional manifold admitting
a WK-spinor. Then $N^3$ is either flat or a space of constant positive
curvature (see [KimF], Theorem 8.2.). Finally, we consider that the case of
$K=M=1$, i.e., $\omega_{12}= \sigma^3$ and $\omega_{23} = \sigma^1$.
Differentiating the equation $\omega_{12}= \sigma^3$, we obtain
$$ \omega_{13} \wedge \omega_{32} - \frac{B}{2} \sigma^1 \wedge \sigma^2 = 
d \omega_{12} = d \omega^3 = \omega_{31} \wedge \sigma^1+
\omega_{32} \wedge \sigma^2 $$
$$ - \frac{B}{2} \sigma^1 \wedge \sigma^2 = - \sigma^1 \wedge \sigma^2 . $$

Consequently, $B=2$ and the tensors $T$ and $Ric$ are given by the 
matrices
$$ T= \left( \begin{array}{ccc} 2-C & 0&0\\
0& 2(C-2) &0\\
0&0& 2-C \end{array} \right) \quad , \quad
Ric= \left( \begin{array}{ccc}
C & 0&0\\
0&2&0\\
0&0&C \end{array} \right) . $$

The second condition of Theorem 1 yields the equations $(S=2+2C)$:
\begin{eqnarray*}
 8 \lambda^2 ( SC - 2C^2 ) - 4 \lambda S(2-C) - S^2 C&=&0 \\
 8 \lambda^2 ( 2S - 8 ) + 8 \lambda S(2-C) - 2 S^2 &=&0 
\end{eqnarray*}

Solving these equations with respect to $\lambda$ and $C$ we 
obtain the three solutions: 

\begin{enumerate}
\item $C=2$ and $\lambda= \pm \frac{3}{2}$. Then $N^3$ is isometric
to $S^3$.
\item $C=-1$ and $\lambda =0$. Then the scalar curvature $S=0$ is zero.
\item $C= \frac{1}{2} (-1 \pm \sqrt{5})$ and $\lambda = 1 \pm \frac{\sqrt{5}}{2}$.
These metrics are the non-Einstein Sasakian metrics on $S^3$ admitting 
WK-spinors (see [KimF]). The corresponding space is contained
in the family $N^3(K,L,M)$.
\end{enumerate}

We have discussed all possibilities and, therefore, we have finished the proof of the 
main Theorem.\\
\mbox{} \hfill \rule{3mm}{3mm}\\

\section{The moduli space of solutions}

The moduli space of all 3-dimensional Riemannian manifolds with constant
scalar curvature $S \not= 0$ and WK-spinors is given by the triples
$\{ K,L,M \}$ of real numbers satisfying the equation of order six 
$Q(K,L,M)=0$. The polynomial $Q$ is symmetric in $\{K, -L, M\}$. Denote by
$$ \gamma_1 =K-L+M \quad , \quad \gamma_2 = - KL+KM-LM \quad , \quad
\gamma_3 = - KLM $$

the elementary symmetric functions of these variables. Then we have
$$ Q= 4 \gamma_1 \gamma_2 \gamma_3 - \gamma_2^3 - 4 \gamma_3^2 . $$

Consider the projective variety $V_{\Bbb C} \subset {\Bbb P}^2 ({\Bbb C})$ defined
by the homogeneous polynomial $Q$:
\[ V_{\Bbb C} = \Big\{ [K:L:M] \in {\Bbb P}^2 ({\Bbb C}) :Q(K,L,M)=0 \Big\} . \]

$V_{\Bbb C}$ has three singular points:
$$ V_{\Bbb C}^{\mathrm{sing}} = \Big\{ [1:0:0] ,\, \, [0:1:0] , \, \,  [0:0:1] \Big\} $$

and these points correspond to the flat metric. We will now parametrize the variety $V_{\Bbb C}$ by two meromorphic functions defined on a smooth Riemann surface. $V_{\Bbb C}$ is given by the equation $(K=1)$:
$$ Q(1,L,M)=L^3(M-1)^2(M+1)+L^2M(1+M)^2-LM^2(1+M)-M^3=0  . $$

Let us introduce the variables
$$ a= M-L-LM \quad , \quad b= (L-M)LM . $$

Then we obtain $Q(1,L,M)=-a^3+4b(1+a)$ and the equation defining the variety $V_{\Bbb C}$ becomes much simpler:
$$ b= \frac{1}{4} \, \, \frac{a^3}{1+a} . $$

Next we consider a square root of $a+1$ and we solve the equations
$$ z^2-1=a=M-L-LM \quad , \quad \frac{1}{4} \, \, \frac{(z^2-1)^3}{z^2} =b=(L-M)LM $$

with respect to $L$ and $M$. Then we obtain four solution pairs $\{L,M\}$ depending on the variable $z$. For example,
$$ L(z)= \frac{-(1+z)(1-2z+z^2+ \sqrt{(1+z)(1+3z-5z^2+z^3}))}{4z} $$
$$ M(z)= \frac{(1+z)(1-2z+z^2+ \sqrt{(1+z)(1+3z-5z^2+z^3}))}{4z}  . $$

The polynomial
$$ (z+1)(1+3z-5z^2+z^3)=(z+1)(z-1)(z+(2+ \sqrt{5}))(z+ (2- \sqrt{5})) $$

has four different zeros. The square root $\sqrt{(1+z)(1+3z-5z^2+z^3)}$ is a meromorphic function on the compact Riemann surface of genus $g=1$. Consequently, there exists a torus ${\Bbb C} / \Gamma$ and elliptic functions $L,M: {\Bbb C}/ \Gamma \to {\Bbb P}^1 ({\Bbb C})$ such that the components of the variety $V_{\Bbb C} \backslash V_{\Bbb C}^{\mathrm{sing}}$ are parametrized by $L$ and $M$:
$$ V_{\Bbb C} = \Big\{ [1: L(z):M(z)]: \, \, \, z \in {\Bbb C} / \Gamma \Big\} . $$

The functions $L-M$ and $L \cdot M$ are given by the formulas:
$$ L-M = - \frac{(1+z)(z-1)^2}{2z} \quad , \quad L \cdot M = - \frac{(1+z)^2(z-1)}{2z} $$

The moduli space we are interested in coincides with the real points 
of the projective variety $V_{\Bbb C}$. If
$K=0$, the only solutions of the equation $Q(0,L,M)=0$ are $L=0$ or $M=0$, i.e.,
the points $[0:1:0]$ and $[0:0:1]$. Therefore we can parametrize the moduli space by
the parameter $M \in {\Bbb R}$ solving the equation $Q(1,L,M)=0$ with respect to
$L= {L} (M)$. In this way we obtain a configuration of six curves in ${\Bbb P}^2({\Bbb R})$ 
connecting the three singular points of $V_{\Bbb C}$ (see the figure in the
Introduction). However,  we obtain geometrically different metrics on $S^3$ only
for two curves parametrized by the real parameter $0 \le M \le \infty$. The graphs
of the function ${L}_{\pm} (M)$ are given in the following figure:\\

\begin{center}

\[
\epsfig{figure=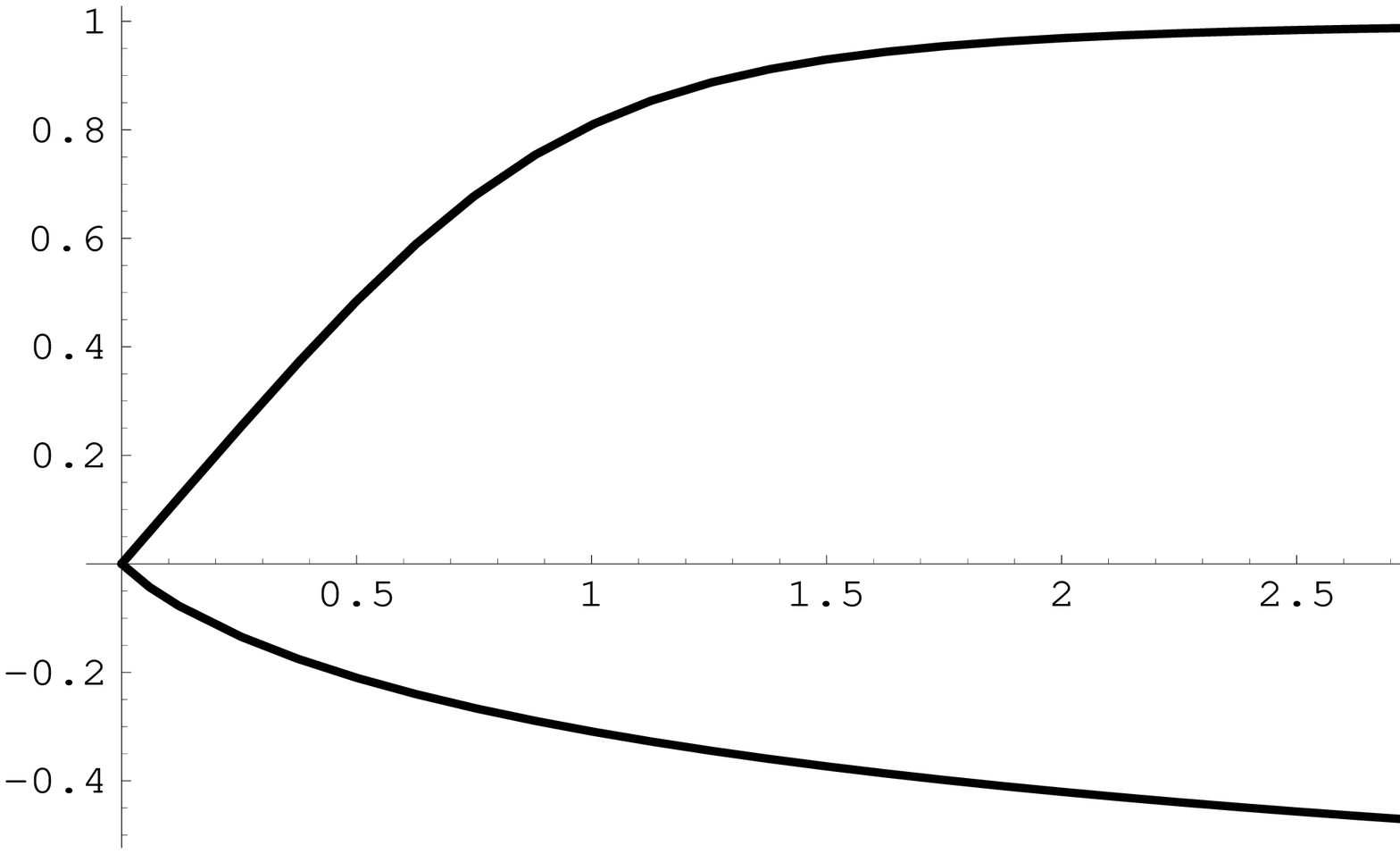,width=8cm}
\]

{\small (Figure 1)}
\end{center}

The functions $L_{\pm} (M)$ are monotone and tend to $\pm 1$ in case that $M$ tends to
infinity. Let us discuss the geometric invariants of these metrics. The graph of the
scalar curvatures $S_{\pm} (M)$ depending on $M$ is given by the next figure:\\

\begin{center}

\[
\epsfig{figure=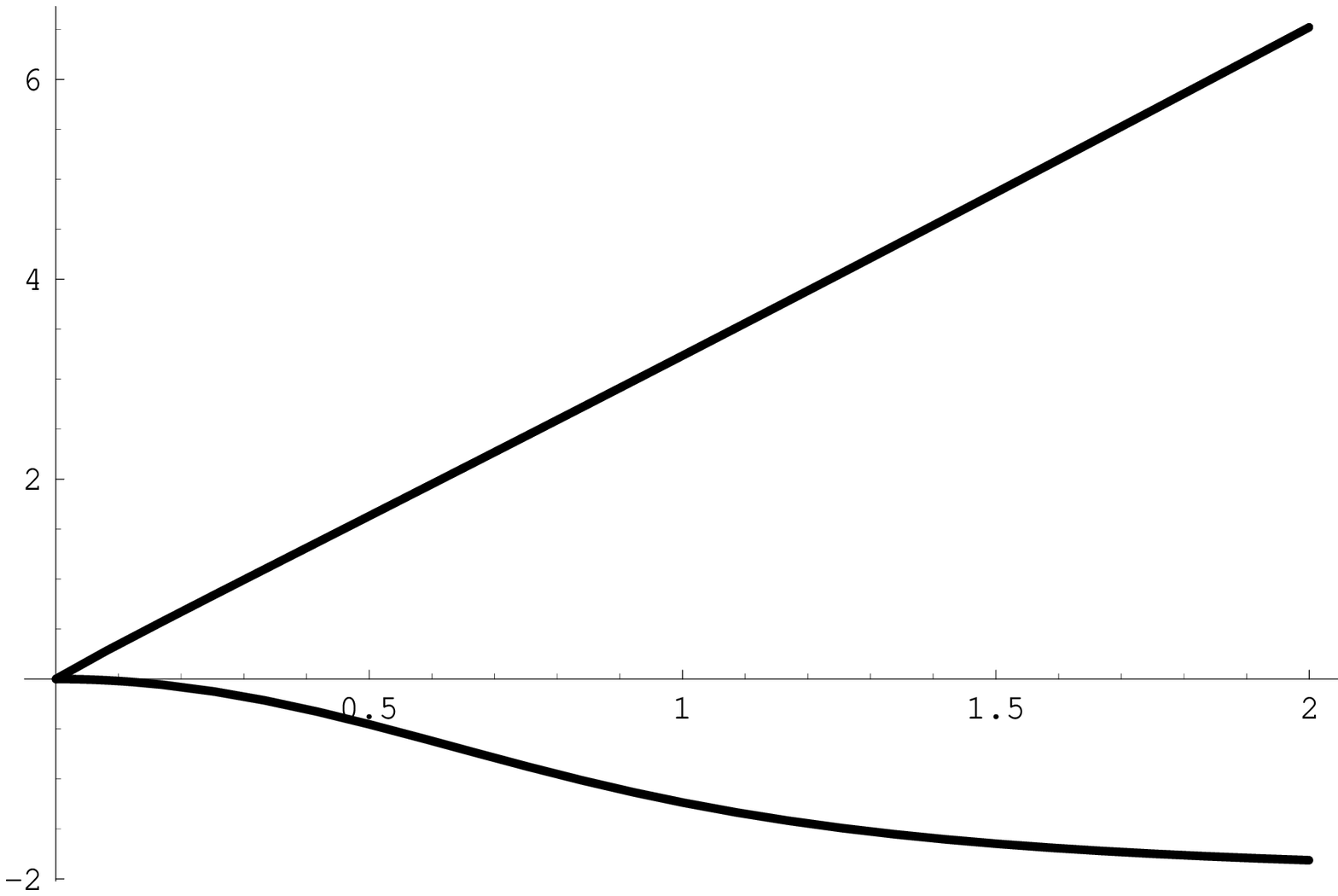,width=8cm}
\]

{\small (Figure 2: The scalar curvatures)}
\end{center}

Next we plot the eigenvalues $A_{\pm} (M), B_{\pm} (M), C_{\pm} (M)$ of the Ricci tensor for both fami\-lies of metrics:

\begin{center}

\[
\epsfig{figure=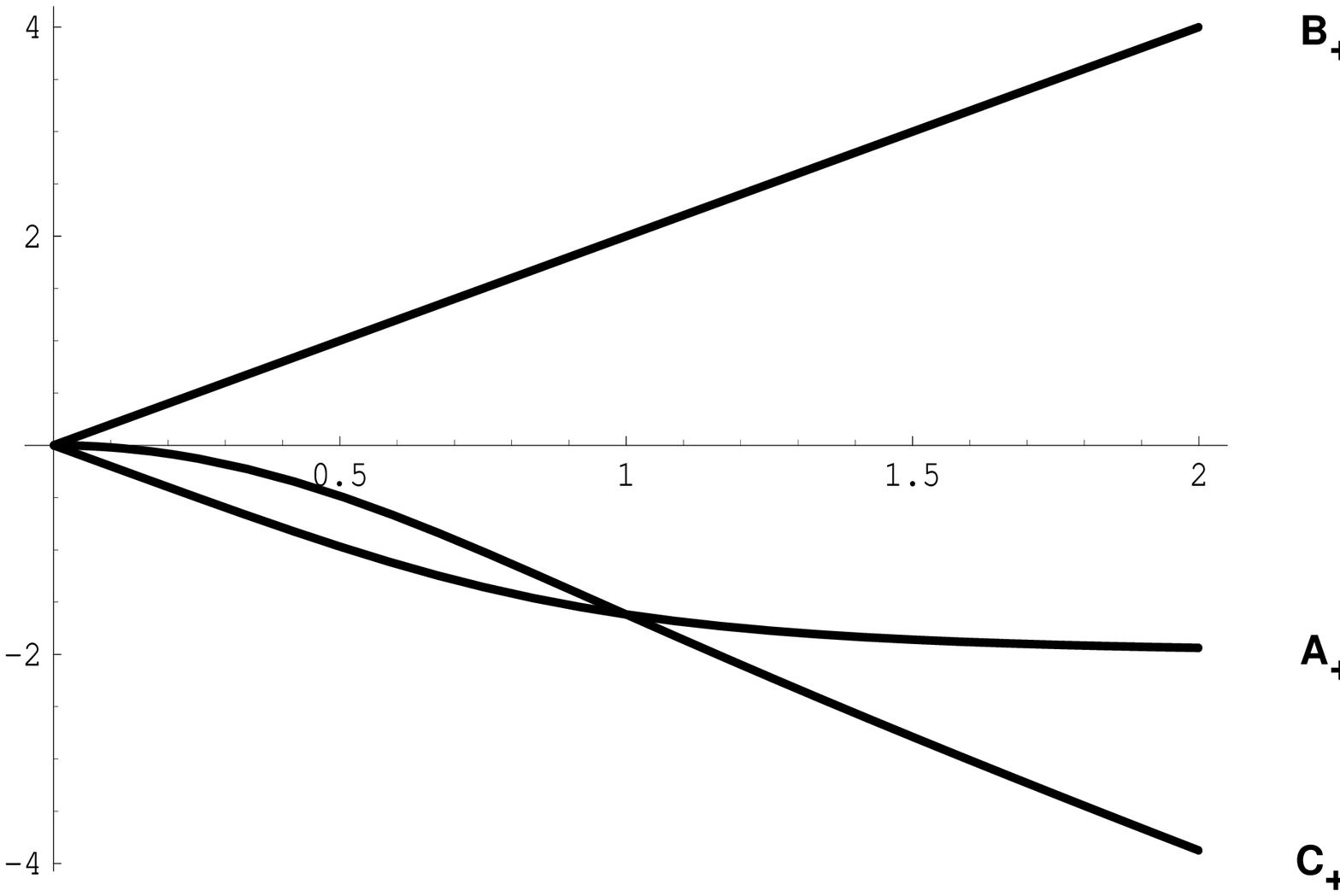,width=8cm}
\]

{\small (Figure 3: The eigenvalues of the Ricci tensor for $L_+(M)$)}
\end{center}

\begin{center}

\[
\epsfig{figure=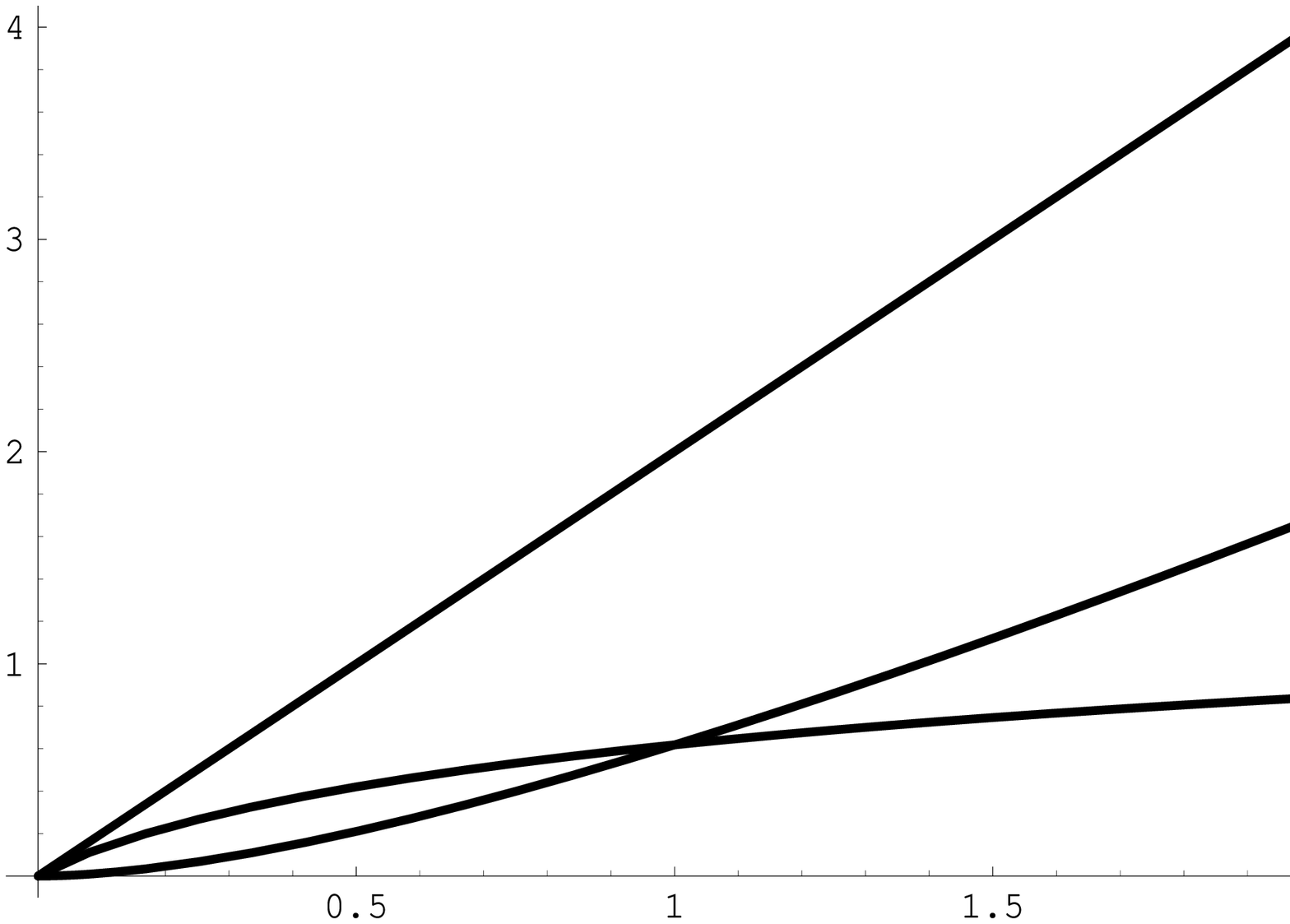,width=8cm}
\]

{\small (Figure 4: The eigenvalues of the Ricci tensor for  $L_- (M)$)}
\end{center}

In dimension $n=3$ the number
$$\lambda^2 (D) \cdot [ \mbox{vol} \, (N^3)]^{\frac{2}{3}} $$ 

is a homothety invariant, where $\lambda (D)$ is an eigenvalue of the Dirac operator. 
In case of a WK-spinor we have
$$ \lambda^2 = \frac{1}{8} \, \, \frac{S^3}{S^2 - 2 |Ric|^2} $$

and, therefore, we obtain the formula
$$\lambda^2 \cdot \mbox{vol}^{\frac{2}{3}} = \frac{1}{8} (2 \pi^2)^{\frac{2}{3}} \frac{S^3}{S^2 - 2 |Ric|^2} \, \, \frac{1}{\{ |K-L||M-L||K+M|\}^{\frac{2}{3}}} . $$

The next figures contain the graph of $\lambda^2 \mbox{vol}^{\frac{2}{3}} (M)$ depending on the parameter $M$ for both families of metrics. \\

\begin{center}

\[
\epsfig{figure=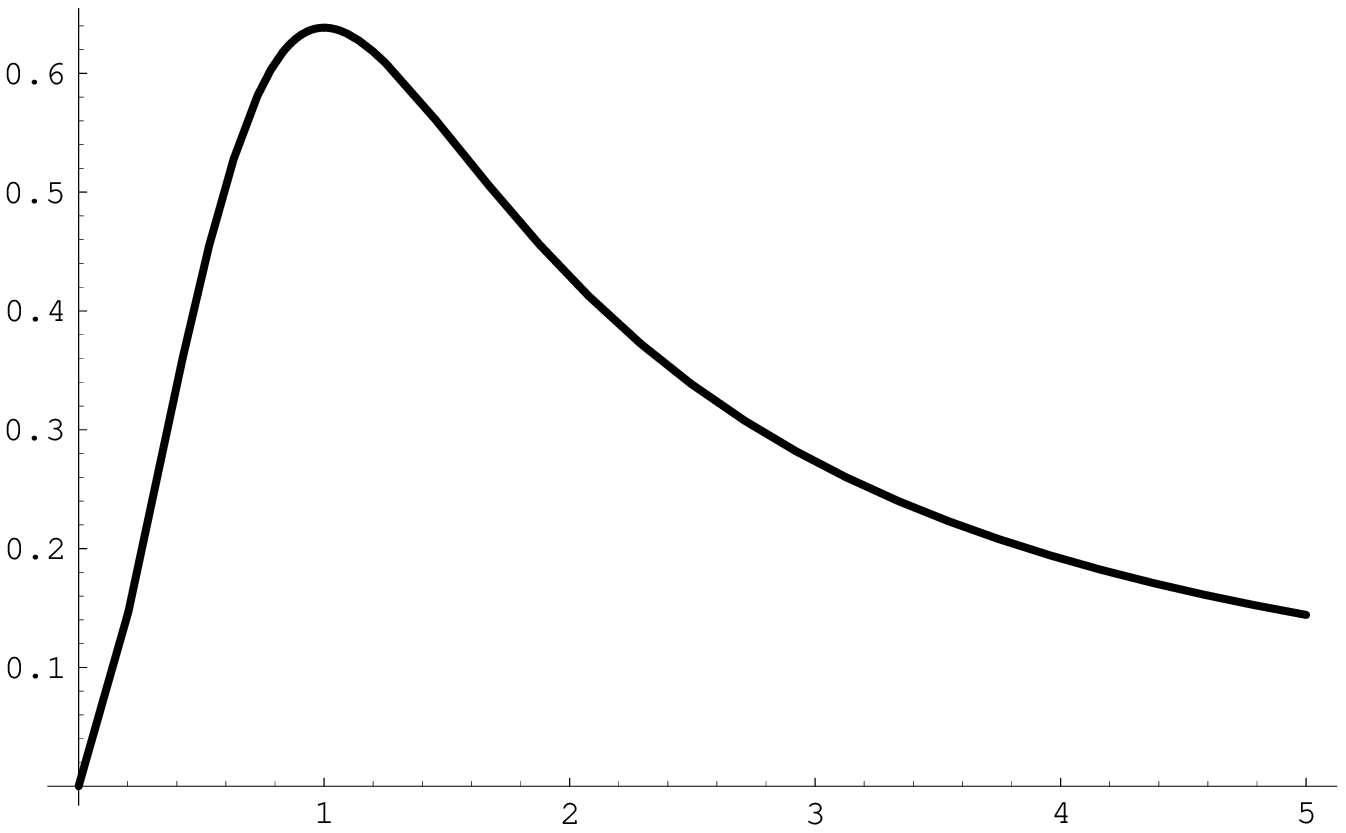,width=8cm}
\]

{\small (Figure 5: $\lambda^2 \mbox{vol}^{\frac{2}{3}}$ in case of $L_+(M)$)}
\end{center}

\begin{center}

\[
\epsfig{figure=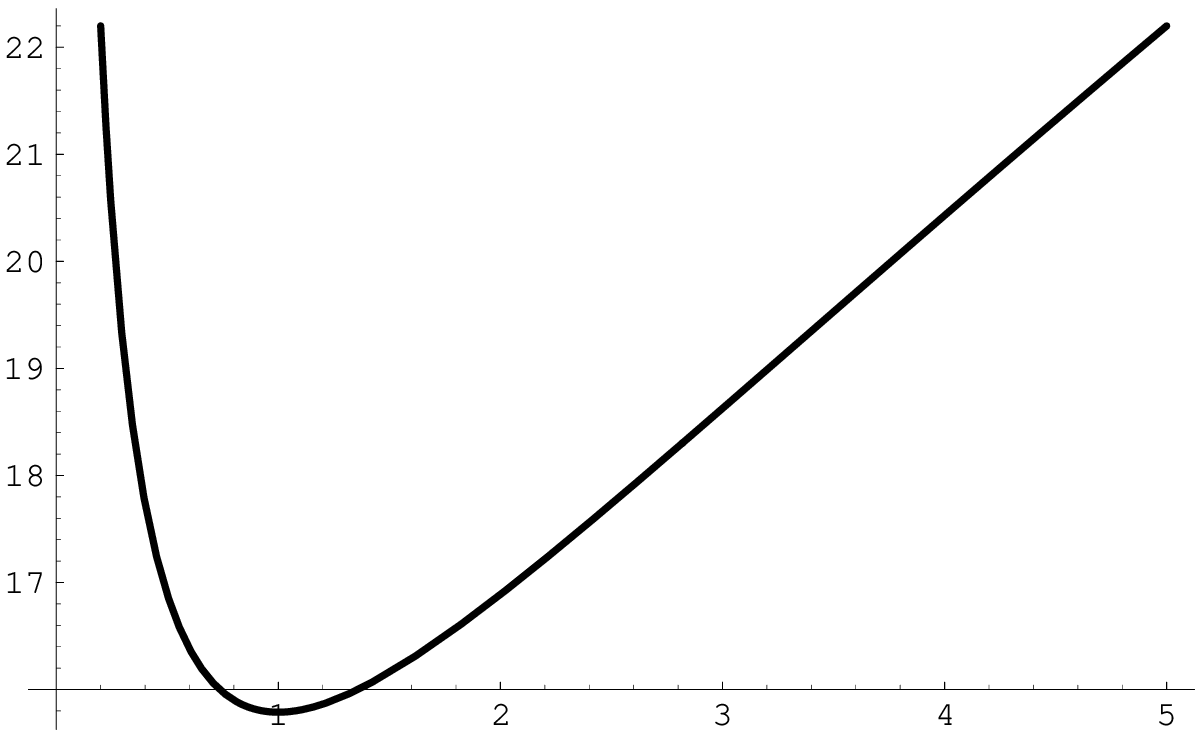,width=8cm}
\]

{\small (Figure 6: $\lambda^2 \mbox{vol}^{\frac{2}{3}}$ in case of $L_-(M)$)}
\end{center}

Finally, let us discuss the behaviour of the rational function
$$ \Psi = \frac{L^2}{KM} $$

on the variety $V_{\Bbb C} \subset {\Bbb P}^2 ({\Bbb C})$. It turns out that $\Psi$ has simple zeros at the singular points $[1:0:0]$ and $[0:0:1]$. Indeed, solving the equation defining $V_{\Bbb C}$ with respect to $L=L(M)$ $(K=1)$ we obtain
$$ \lim\limits_{M \to 0} \frac{L^2(M)}{M} =0 \quad , \quad \lim\limits_{M \to 0} \frac{d}{dM} \Big( \frac{L^2(M)}{M} \Big) = 1 . $$

The third singular point $[0:1:0]$  is a pole of order two. In the regular part of $V_{\Bbb C}$ the function $\Psi$ has 12 ramification points. Among them 10 points are first order ramification points. The ramification points of order two are the points
$$ [K:L:M]=[1: \frac{1}{4} (1 \pm \sqrt{5}):1 ] . $$

These parameters correspond precisely to the non-Einstein Sasakian metrics on $S^3$ admitting solutions of the Einstein-Dirac equation.\\


\bigskip


\begin{thebibliography}{MMMM}
\bibitem[F] {} Th. Friedrich, Dirac-Operatoren in der Riemannschen Geometrie,
Vieweg Verlag Wiesbaden 1997.
\bibitem[FKim]{} Th. Friedrich and E.C. Kim, Some remarks on the Hijazi inequality and 
gene\-ralizations of the Killing equation for spinors, to appear in Journ. Geom. Phys. (2000).
\bibitem[KimF] {} E.C. Kim and Th. Friedrich, The Einstein-Dirac equation on Riemannian 
spin manifolds, Journ. Geom. Phys. 33 (2000), 128-172.
\end{thebibliography}
\end{document}